\documentclass[12pt]{article}
\usepackage{amsmath,amssymb,latexsym,graphicx}
\newtheorem{theorem}{Theorem}
\newtheorem{lemma}[theorem]{Lemma}
\newtheorem{proposition}[theorem]{Proposition}
\newtheorem{remark}[theorem]{Remark}

\newtheorem{definition}[theorem]{Definition}

\def\pf{{\bf Proof }}%^%

\begin{document}
\title{On polynomially integrable domains in Euclidean spaces}
%\address{Department of Mathematics, Bar-Ilan University, Ramat-Gan, 52900; Department of Mathematics, Holon Istitute of Technology, Holon; Israel }
%\acknowledgement {The author thanks Mikhail Zaidenberg for drawing the author's attention to the subject of this article and stimulating initial discussions.}
%The work of the second author was supported in part by the NSF Grants DMS 9971674 and
%0002195.}
\author{Mark~ L.~Agranovsky}
\maketitle
%\begin{center}
%Bar-Ilan University, Holon Institute of Technology
%\end{center}
%\begin{\dedication}

%\bigskip
%\hspace{5cm}
%{\it To the memory of my  friend Sasha Vasiliev}
%\bigskip

\begin{abstract}
Let $D$ be a bounded domain in $\mathbb R^n,$ with smooth boundary. Denote $V_D(\omega,t), \ \omega \in S^{n-1}, t \in \mathbb R,$ the Radon transform of the characteristic function $\chi_{D}$ of the domain $D,$ i.e.,
the $(n-1)-$ dimensional volume of the intersection $D$ with the hyperplane $\{x \in \mathbb R^n: <\omega,x>=t \}.$
If the domain $D$ is an ellipsoid, then the function $V_D$ is algebraic and if, in addition, the dimension $n$ is odd,  then $V(\omega,t)$ is a polynomial with respect to  $t.$  Whether odd-dimensional ellipsoids are the only bounded smooth domains with  such a property? The article is devoted to partial verification and discussion of that question.
\end{abstract}

\section{Introduction}\label{S:Intro}

V. A. Vassiliev proved \cite{Vas1,Vas} that  if $D$ is a bounded domain in even-dimensional space $\mathbb R^{n},$ with $C^{\infty}$ boundary, then the two-valued function,  evaluated the $n-$ dimensional volumes $\widehat V^{\pm}(\omega,t)$ of the two complementary portions of $D$ which are cut-off by the section of $D$ by the hyperplane $\{<\omega,x>=t\},$ is not an algebraic function of the parameters of the hyperplane.
That means that  no functional
equation $Q(\omega_1,\cdots,t,\omega_n, \widehat V^{\pm}(\omega,t))=0,$ where $Q$ is a polynomial, is fulfilled. The result of \cite{Vas} is a multidimensional generalization of a celebrated Newton's Lemma XXVIII in Principia \cite{N}  for convex ovals in $\mathbb R^2.$ The smoothness condition is essential.

The main object of interest in this article is the {\it section-volume function}
$$V_D(\omega,t)=\int\limits_{<x,\omega>=t}f(x)dV_{n-1}(x),$$
which is the $t-$ derivative of $\widehat V^{\pm}.$

Contrary to the multi-valued cut-off-volume function $\widehat V_D,$ its derivative $V_D=\frac{d}{dt}\widehat V_D$ is  single-valued and can be algebraic in any dimension of the ambient space
(the differentiation can "improve" the algebraicity, just keep in mind the algebraic function $\frac{1}{1+x^2}$ which is the derivative of the non-algebraic function $\tan x.$)

For example, if $D$ is the unit ball in $\mathbb R^{n}$ then
$$V_D(\omega,t)=c(1-t^2)^{\frac{n-1}{2}}$$
is algebraic in any dimension. If $n$ is odd then, even better, $V_D$ is a polynomial in $t.$ Here and further, $c$ will denote a nonzero constant, which exact value is irrelevant.

Ellipsoids in odd-dimensional spaces  also have the same property. For instance, if
$$E=\{\sum\limits_{j=1}^n \frac{x_j^2}{b_j^2}=1\}$$
then
$$V(\omega,t)= c \frac{b_1 \cdots b_n}{h^n(\omega)} (h^2(\omega)- t^2)^{\frac{n-1}{2}},$$
where
$$h(\omega)=\sqrt{\sum\limits_{j=1}^n b_j^2\omega_j^2}.$$

\begin{definition}
We will call a domain $D$ {\it polynomially integrable} if its section-volume function coincides with a polynomial in $t:$
 $$V_D(\omega,t)=\sum\limits_{j-1}^N a_j(\omega) t^j$$
in the domain $V_D(\omega,t) > 0$ of all $(\omega,t)$ such that
the hyperplane $<\omega, x>=t$ hits the domain $D.$
\end{definition}
We assume that the leading coefficient $a_N(\omega)$ is not identical zero  and in this case we call $N$ {\it the degree} of the polynomially integrable domain $D.$

Note that in \cite{Vas} the term "algebraically integrable" refers to the cut-off volume function $\widehat V_D,$ rather than to $V_D$
but polynomial integrability with respect to either function is the same.

It was conjectured in \cite{Vas} that  ellipsoids are the only smoothly bounded algebraically  integrable domains in odd-dimensional Euclidean spaces. We are concerned with the similar question about
the section-volume function $V_D(\omega,t)$ and the condition is that this function is a polynomial in $t.$ Although polynomials are definitely algebraic functions, we do not impose conditions on the dependence on $\omega.$   Still, our conjecture is that the odd-dimensional ellipsoids are the only bounded polynomially integrable domains in odd-dimensional Euclidean spaces.

The present article contains some partial results and observations on the conjecture.

Namely, in Sections 1-4 we prove that there are no bounded polynomially integrable domains with smooth boundaries in
Euclidean spaces of even dimensions (Theorem \ref{T:even}), while in odd dimensions the following is true: 2) polynomially integrable bounded domains with smooth boundaries are convex (Theorem \ref{T:convex}),
3) there are no polynomially integrable bounded domains in $\mathbb R^n,$ with smooth boundaries, of degree strictly less than $n-1$ (Theorem \ref{T:degree}),
4) polynomially integrable bounded domains in $\mathbb R^n,$ with smooth boundaries,  of degree $\leq n-1$ are only ellipsoids (Theorem \ref{T:T1}, 4)
polynomially integrable bounded domains in $\mathbb R^3,$ with smooth boundaries, having axial and central symmetry are only ellipsoids (Theorem \ref{T:axial}).
In Section \ref{S:Stat} a relation of polynomial integrability and stationary phase expansions is discussed. Some open questions are formulated in Section \ref{S:Open}.

\section{Preliminaries} \label{S:Prelim}
\subsection{Support functions}

For any bounded domain $D$ define the support functions

\begin{equation}\label{E:h}
\begin{aligned}
h^{+}_D(\omega)&=\sup\limits_{x \in D} <\omega,x>, \\
h^{-}_D(\omega)&=\inf\limits_{x \in D} <\omega,x>,
\end{aligned}
\end{equation}
where $\omega$ is a  vector in the unit sphere $S^{n-1}$ and $\langle, \rangle$ is the standard inner product in $\mathbb R^n$

The two support functions are related by
$$h_{D}^{+}(-\omega)=-h_{D}^{-}(\omega).$$
The relation with translations is given by
\begin{equation}\label{E:D+a}
h_{D+a}^{\pm}(\omega)=h_{D}^{\pm}(\omega)+<a,\omega>.
\end{equation}

 If $ t \notin [h_D^{-}(\omega),h_D^+(\omega)]$ then the  $<\omega,x>=t$ is disjoint from the closed domain $\overline D$ (and therefore $V_D(\omega,t)=0$ ).

\includegraphics[width=14cm]{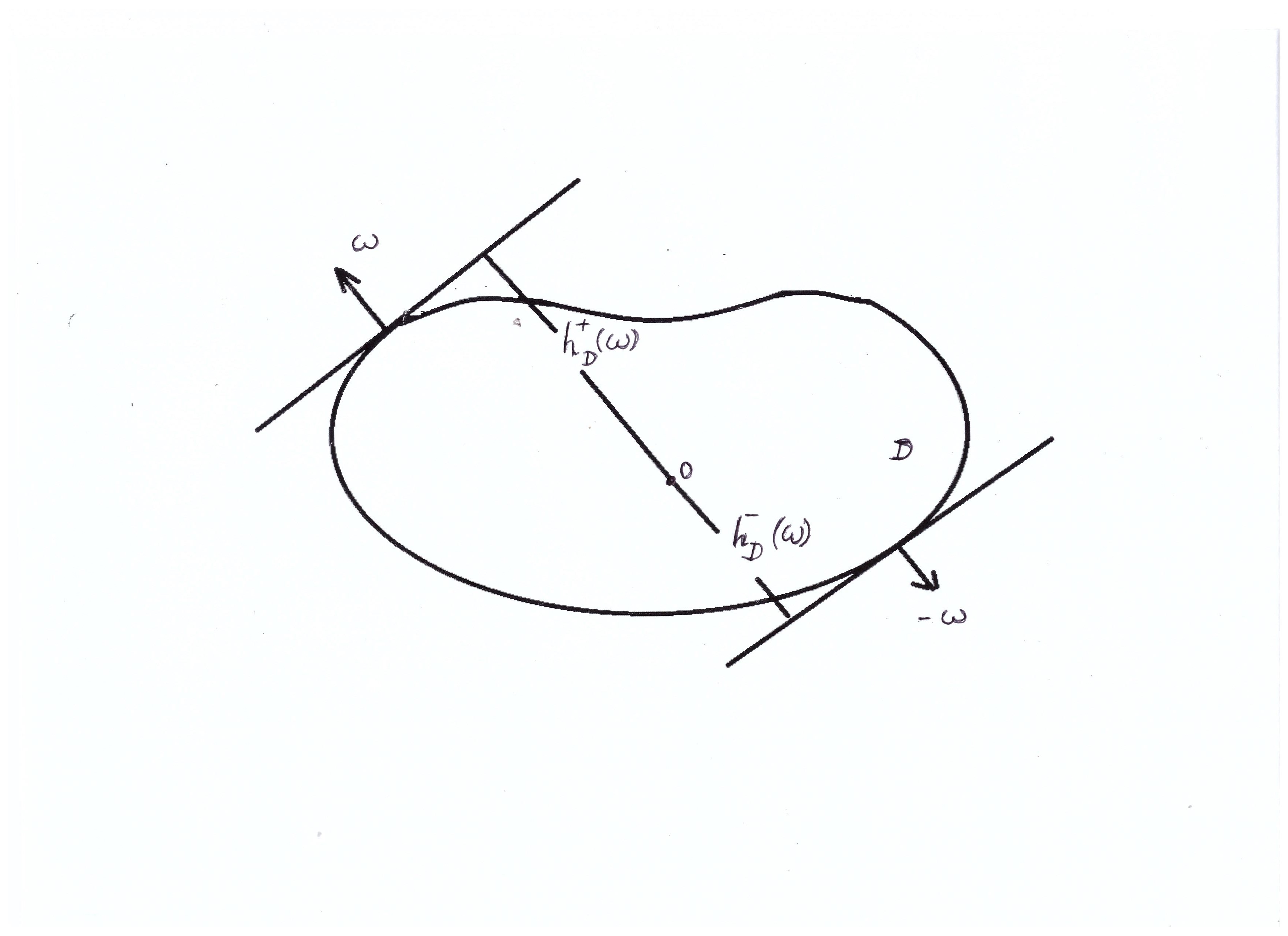}

\subsection{Radon transform}
We will be using several standard facts about the Radon transform. They can be found in many books or articles on Radon transform, see e.g.,
\cite{GG},\cite{Helg}.

The Radon transform of a continuous compactly supported function $f$ in $\mathbb R^n$ is defined as the integral
$$Rf(\omega,t)=\int\limits_{<\omega,x>=t}f(x)dV_{n-1}(x)$$
of $f$ over the $(n-1)-$ dimensional plane $<\omega,x>=t,$ with the unit normal vector $\omega,$ on the distance $t$ from the origin, against the $(n-1)-$ dimensional volume element
$dV_{n-1}.$

There is a nice relation, called {\bf Projection-Slice Formula}, between  Fourier and Radon transforms:
\begin{equation}\label{E:SPF}
\widehat f(r\omega)=F_{t \to r} Rf(\omega,t),
\end{equation}
where $\widehat f$ is the $n-$ dimensional Fourier transform and $F_{t \to r}$ is the one-dimensional Fourier transform with respect to the variable $t.$
The formula follows immediately by writing the integral as a double integral against $dV_{n-1}$ and $dt.$

The Radon transform can be inverted in different ways. The following formula is called back-projection inversion formula (\cite{Helg}, Thm 2.13) :
$$f(x)=c \Delta^{\frac{n-1}{2}} \int\limits_{\omega \in S^{n-1}} Rf(\omega,<\omega,x>) d\omega,$$
where $f$ is  sufficiently smooth compactly supported function.

There is a crucial difference between the inversion formulas in even and odd dimensions. If $n$ is even then the exponent $\frac{n-1}{2}$ is fractional and $\Delta^{\frac{n-1}{2}}$
becomes an integral operator, so that the inversion formula is not local. However, when $n$ is odd then $\frac{n-1}{2}$ is integer and the operator in front of the integral  is differential so that
the formula is local.

Due to the relation
\begin{equation}\label{E:diff}
\widehat {\frac{\partial f}{\partial x_i}}(\omega,t)=\omega_i \frac{\partial \widehat f}{dt}(\omega,t),
\end{equation}
{\bf the inversion formula}  can be rewritten in odd dimensions as:
\begin{equation}\label{E:inversion}
f(x)=c \int\limits_{\omega \in S^{n-1}} \left(\frac{d^{n-1}}{dt^{n-1}} Rf \right)(\omega, <\omega,x>) d\omega,
\end{equation}
where $d\omega$ is the Lebesgue measure on the unit sphere $S^{n-1}.$
The Plancherel formula for the Radon transform in odd dimensions is:
\begin{equation}\label{E:Planch}
c \int\limits_{\mathbb R^n} f(x)g(x)dV(x)=\int\limits_{S^{n-1} \times \mathbb R} \left(\frac{ d^{n-1} }{dt^{n-1}} Rf\right)(\omega,t)Rg(\omega,t)d\omega dt.
\end{equation}
The following conditions characterize the range of Radon transform in class of Schwartz functions:
\begin{enumerate}
\item
$g(-\omega,-t)=g(\omega,t),$
\item The $k-$  moment
$\int\limits_{\mathbb R^n}g(\omega,t)t^kdt, k =0,1,\cdots$
extends from the unit sphere $S^{n-1}$ to $\mathbb R^n$ as a polynomial of degree at most $k.$
\end{enumerate}
The immediate corollary of property 1 is that if $V_D(\omega,t)=\sum\limits_{k=0}^N a_k(\omega)t^k$ then
$$a_k(-\omega)=(-1)^k a_k(\omega),$$
i.e. $a_k$ is an even function on $S^{n-1}$ when $k$ is even and is an odd function when $k$ is odd.

In the sequel, the above facts will be used for the section-volume function $V_D$ which is just the Radon transform of the characteristic function of the domain $D:$
$$V_D(\omega,t)=(R\chi_D)(\omega,t), \ t \in [h_D^{-}(\omega),h_D^+(\omega)].$$

All domains under consideration will be assumed bounded, with $C^{\infty}$ boundary, although some statements are true even under weaker smooth assumptions.
As it is shown in \cite{Vas},\cite {Vas1}, smoothness plays an important role in that circle of questions.

\section{There are no polynomially integrable domains in even dimensions}
The condition of algebraic integrability in \cite{Vas} involves both variables $\omega$ and $t,$ while the polynomial integrability imposes a condition with respect to $t$ only, so that  the statement formulated in the title on this section is not a straightforward corollary  of the result in \cite{Vas}. However, this statement easily follows from the asymptotic behavior of the section-volume function near the boundary points of $D:$
\begin{theorem}\label{T:even} There are  no polynomially integrable domain with $C^{2}-$ smooth boundary  in $\mathbb R^n$ with even $n.$
\end{theorem}
\pf Let $D$ be such a domain. We will show that $V_D(\omega,t)$ can not behave polynomially when the section of $D$ by the hyperplanes $<x,\omega>=t$ shrink to an elliptical point on $\partial D.$ The elliptical point of the surface $\partial D$ is  a point at which the principal curvatures are all nonzero and of the same sign.

To find such a point, consider  the maximally distant, from the origin, boundary point $a \in \partial D:$
$$|a|=max_{x \in \partial D} |x|.$$ Using rotation and translation, we can move $a$ to $0$ and make the tangent plane at 0 the coordinate plane $$T_{\partial D}(0)=\{x_n=0\}$$
so that the domain $D$ is in the upper half-plane $x_n \geq 0.$  Let $b$ the image of the point $a$ under those transformations
Then $b$ is located on the $x_n-$ axis and
$$D \subset B(b,|b|),$$
where $B(b,r)$ is the ball of radius $r$ with center $b.$
The boundary surface $\partial D$ can be represented, in a neighborhood of  $0 \in \partial D,$ as the graph
$$\partial D \cap U=\{x_n =f(x^{\prime}), x^{\prime}=(x_1,...,x_{n-1})\}.$$
Since the tangent plane at $0$ is $x_n=0$ the first differential vanishes at the origin, $df_0=0.$ Using rotations around the $x_n$ axis we can diagonalize the second differential and then the equation of $\partial D$ near $0$ becomes
$$\partial D \cap U=\{x_n= \frac{1}{2} \sum\limits_{j=1}^{n-1}\lambda_j x_j^2 + S(x^{\prime}) \},$$
where $S(x^{\prime})=o(|x^{\prime}|^2), x^{\prime} \to 0.$

The coefficients $\lambda_j$ are the principal curvatures of the surface $\partial D$ and since $D \subset B(b,|b|)$ all of them are not less
than the curvature $\frac{1}{|b|}$ of the sphere $S(b,|b|):$
$$ \lambda_j  \geq \frac{1}{|b|} >0, j=1,...,n-1.$$
Indeed,
the equation of the sphere $S(b,|b|)$ is
$$x_n=|b|-\sqrt{|b|^2-|x^{\prime}|^2}$$
and since $D \subset B(b,|b|)$ we have
$$
\begin{aligned}
\frac{1}{2} \sum\limits_{j=1}^{n-1}\lambda_j x_j^2 +o(|x^{\prime}|^2) \geq |b|-\sqrt{|b|^2-|x^{\prime}|^2}&= \\
&=\frac{|x^{\prime}|^2}{|b|+\sqrt{|b|^2-|x^{\prime}|^2}}.
\end{aligned}
$$
and then the inequalities follow by the passage to the limit.

 Let us turn further to the function $V_D.$ The value $V_D(e_n,t), e_n=(0,...,0,1), t>0,$ is the $(n-1)-$ dimensional volume of the section $D \cap \{x_n=t\}$ and the leading term of the asymptotic when $t \to 0$ is defined by the leading, quadratic, term in the decomposition of $f.$ This leading term of $V_D(e_n,t)$ is the volume of the ellipsoid
$$\frac{1}{2} \sum\limits_{j=1}^{n-1}\lambda_j x_j^2=t,$$
which is $c \frac{t^{\frac{n-1}{2}}}{\lambda_1 \cdots \lambda_{n-1}}, c= \frac{(2\pi)^{\frac{n-1}{2}}}{\Gamma(\frac{n+1}{2})}.$
Therefore,
$$V_D(e_n,t)=c \frac{t^{\frac{n-1}{2}}}{\lambda_1 \cdots \lambda_{n-1}}+ o(t^{\frac{n-1}{2}}), t \to 0+.$$
Thus, the function $V_D(e_n,t)$ has zero at $t=0$ of  order $\frac{n-1}{2}.$ If this function is a polynomial then the order must be integer, which is not the case when $n$ is even.
Thus, $V_D(e_n,t)$ is not a polynomial in $t.$ In fact, all we have used is that $a$ is an elliptical point of the convex hypersurface $\partial D.$ Since the points of $\partial D,$ sufficiently close to $a,$
are elliptical as well, we conclude that, moreover, for an open set of directions $\omega,$ the function $V_D(\omega,t)$ is not a polynomial in $t.$  This proves the theorem.

\section{Polynomially integrable domains in $\mathbb R^{2m+1}$ are convex}

\begin{lemma} \label{L:a_k}
Suppose that the section-volume function $V_D$ coincides with a polynomial in $t:$
$$V_D(\omega,t)=a_0(\omega)+a_1(\omega)t +...+ a_N(\omega) t^N,$$
when $t \in [h_D^-(\omega),h_D^+(\omega)].$
Then
$$\int\limits_{|\omega|=1} a_{k}(\omega)p(\omega) dA(\omega)=0$$
for any polynomial $p$ of $\deg \ p \leq k-n+1$ and any $k>n-1.$
\end{lemma}

\pf Regarding the function $V_D$ as the Radon transform of the characteristic function $\chi_D$ of the domain $D,$ we can write
by the inversion formula (\ref{E:inversion}):
$$
\begin{aligned}
&\chi_D(x)=\int\limits_{|\omega|=1}\frac{d^{n-1}}{dr^{n-1}}V_D(\omega, <x,\omega>) dA(\omega)\\
&=\sum_{k=n-1}^N\frac{k!}{(k-n+1)!} \int\limits_{|\omega|=1}  a_k(\omega)  <x,\omega>^{k-(n-1)}dA(\omega). \\
\end{aligned}
$$

Since $\chi_D(x)=1$ for $x \in D,$ the power series in $x$ in the left hand side equals 1 identically and hence each term with $k-(n-1) >0$ vanishes:
$$\int\limits_{|\omega|=1} a_{k}(\omega)<x,\omega>^{k-(n-1)} dA(\omega)=0, x \in D.$$
Since $x$ is taken from an open set, the functions $y \to <x,y>^{k-n+1}$ span the space of all homogeneous polynomials of degree at most $k-n+1.$
Those polynomials, being restricted on the unit sphere, generate restrictions of all polynomials $p, \ \deg p \leq k-n+1.$ Lemma is proved.

\begin{remark} The assertion of Lemma remains true also under the assumption that $V_D$ expands  as an infinite  power series in $t:$
$$V_D(\omega,t)=\sum\limits_{k=0}^{\infty}a_j(\omega)t^k,$$
in a neighborhood of $[h^{-}(\omega), h^{+}(\omega)].$

\end{remark}
\begin{theorem} \label{T:convex} If a smoothly bounded domain $D$ in $\mathbb R^n,$ with $n$ odd, is polynomially integrable then it is convex.
\end{theorem}

\pf
Let $D$ be such a domain and denote $\widehat D$ the convex hull of  $D.$ If $D \neq \widehat D$  then there is an open portion $\Gamma$ of the boundary $\partial D$ which is disjoint from a neighborhood of $\partial \widehat D.$ Then one can construct a function which vanishes near $\partial \widehat D$ and behaves in a prescribed way on $\Gamma.$ In particular, one can construct a compactly supported $C^{\infty}$ function $\psi$ in $\mathbb R^n$
with the following properties:

\begin{enumerate}
\item
$supp \  \psi \subset \widehat D,$
\item
$\int\limits_{\partial D} \psi(x) \nu_1(x) dS(x) \neq 0,$
where $\nu(x), x \in \partial D, $ is the external unit normal field on  the boundary of $D.$
\end{enumerate}

Denote $\Psi(\omega,t)$ the Radon transform of $\psi:$
$$\Psi(\omega,t)=(R\psi)(\omega, t)=\int\limits_{<\omega,x>=t} \psi(x)dV_{n-1}(x).$$
Formula \ref{E:diff} implies
\begin{equation}\label{E:deriv}
\Psi_1(\omega, t)=\omega_1 \Psi^{\prime}_t (\omega,t).
\end{equation}

Using (\ref{E:deriv}), write the Plancherel formula  (\ref{E:Planch}) for the pair of functions: $\chi_D$ and $\frac{\partial \psi}{\partial x_1}.$
Although formula (\ref{E:Planch}) works for Schwartz functions, we can extend the equality to our case approximating $\chi_D$ by smooth functions, if the differentiation in $t$ in the left hand side  applies to the second, smooth, factor $\Psi=R\psi.$ We obtain

\begin{equation}\label{E:Planch1}
\begin{aligned}
\int_D \frac{\partial \psi}{\partial x_1}(x)dV(x)&=\int\limits_{\mathbb R^n} \chi_D(x)\frac{\partial \psi}{\partial x_1}(x)dV(x)= \\
&=\int\limits_{|\omega|=1}\int_{\mathbb R} R\chi_D(\omega,t) (\Psi_1 (\omega,t))^{(n-1)}_t dA\omega) dt.
\end{aligned}
\end{equation}

By the condition, the Radon transform $R\chi_{D}(\omega, t)$ of the characteristic function of $D$ coincides with the polynomial $V_D(\omega,t)$ as long as the hyperplane $<\omega,x>=t$ hits $D.$
The latter occurs,  if and only if it hits $\widehat D.$ In turn, that is equivalent to the inequality
$h_{-}(\omega) \leq t \leq h_{+}(\omega).$ Thus, the integration in $t$ in the right hand side of (\ref{E:Planch1}) is performed on the segment $[h_D^{-}(\omega), h_D^{+}(\omega)].$

%Altogether, we conclude that  Radon transform $R\chi D$ coincides with the polynomial $V_D(\omega,t)$ in the domain of integration in %(\ref{E:Planch}).
The left hand side of (\ref{E:Planch1})  can be rewritten, using
(\ref{E:deriv}) and Green's formula, as:
$$
\int_D \frac{\partial \psi}{\partial x_1}(x)dV(x)=\int_{\partial D} \psi(x) \nu_1(x) dS(x),
$$
where $\nu_1(x)$ is the first coordinate of $\nu(x)$.

Finally,  (\ref{E:Planch1}) takes the form:
$$
\begin{aligned}
&I:=\int_{\partial D} \psi(x) \nu_1(x) dS(x)=\\
&=\int\limits_{\omega \in S^{n-1}}\int\limits_{h_(\omega)}^{h_+(\omega)}V_D(\omega,t) \omega_1 \Psi^{(n)}(\omega,t)dt dA(\omega).
\end{aligned}
$$

Furthermore, the function $(\psi)^{\prime}_{x_1}$ vanishes in a neighborhood of the boundary $\partial \widehat D.$ Therefore, its integrals over hyperplanes $<\omega,x>=t$ with
$$t \notin I_{\varepsilon}:= [h_{-}(\omega)+ \varepsilon, h_{+}(\omega)-\varepsilon], $$
vanish as well for sufficiently small  $\varepsilon >0$ and
\begin{equation}\label{E:Psi=0}
\Psi (\omega,t)=0, \ \mbox{for} \  t \notin I_{\varepsilon}.
\end{equation}
Hence, all the derivatives of $\Psi(\omega,t)$ vanish at the end points $t =h_{\pm}(\omega)$ of the integration interval and then integration by parts in $t$ yields:
 $$ I=(-1)^n \int\limits_{\omega \in S^{n-1}}\int\limits_{h_(\omega)}^{h_+(\omega)}V_D^{(n)}(\omega,t) \omega_1 \Psi(\omega,t)dt dA(\omega).$$
 More explicitly:
 $$I=(-1)^n \sum_{k=n}^N \frac{k!}{(k-n)!}t^{k-n}\int\limits_{\omega \in S^{n-1}}\int\limits_{-\infty}^{+\infty} a_k(\omega) \omega_1 \Psi(\omega,t)dt d\omega.$$
We have replaced here the integral on $[h_D^{-}(\omega),h_D^{+}(\omega)]$ by integration on $(-\infty,\infty),$ due to \ref{E:Psi=0}.
Consider an arbitrary term in the right hand side:
$$I_k:=\int\limits_{\omega \in S^{n-1}}\int\limits_{\mathbb R} a_k(\omega) t^{k-n} \omega_1  \Psi(\omega,t)dt d\omega,$$
$k=n,n+1,...N.$
The function $\Psi$ is the Radon transform of a smooth compactly supported function $\psi.$ Therefore, according to the description
of the range of Radon transform in Section \ref{S:Prelim}, the moment condition 2) is fulfilled. This means that the function
$$Q_k(\omega):= \int\limits_{\mathbb R}\Psi(\omega,t)t ^{k-n}dt$$
extends from the unit sphere $|\omega|=1$ as a  polynomial of degree at most $k-n.$  Therefore, the function $\omega_1Q_k(\omega)$ extends as a polynomial of degree at most $k-n+1.$
But then
$$ I_k=\int\limits_{\omega \in S^{n-1}}a_k(\omega) \omega_1Q_k(\omega) d\omega=0$$
by Lemma \ref{L:a_k}.

Thus, we have proven that all $I_k=0$ and therefore $I$ which is a sum of $I_k$ is zero as well, $I=0.$ But $I=\int_{\partial D} \psi(x) \nu_1(x) dS(x) \neq 0$ by the choice of the function $\psi.$ This contradiction is obtained under assumption that $D$ is strictly smaller than its convex hull $\widehat D.$ Therefore, the opposite, they coincide and thus $D$ is convex. Theorem is proved.

\begin{remark} It can be seen from the proof that Theorem \ref{T:convex} is valid also in the case when  $V_D (\omega,t)$ decomposes into infinite power series in $t.$
More precisely, if for every fixed $\omega \in S^{n-1}$ the function $V_D(\omega,\cdot)$ extends,  in a neighborhood of the segment $[h_D^{-1}(\omega), h_D^+(\omega)],$
as a power series.
\end{remark}

\section{Polynomially integrable domains in $\mathbb R^{n}, n=2m+1,$ of  degree $\leq n-1$  are ellipsoids }

In this section, we characterize odd-dimensional
ellipsoids by the property of the functions $V_D(\omega,t)$ being polynomials  in $t$ of degree at most $n-1.$ In fact, $n-1$ is also the lower bound for the degrees of those polynomials.
Everywhere in this section the dimension $n$ is assumed being odd.
\begin{theorem} \label{T:degree} If $D$ is a polynomially integrable domain in $\mathbb R^n$ with smooth boundary then for almost all $\omega \in S^{n-1}$ the degree of the polynomial $V_D(\omega, \cdot )$ is at least $n-1.$ In particular, there are no polynomially integrable domains of degree less than $n-1.$
\end{theorem}

\pf
For every point $x \in \partial D$ denote by $\kappa(x)$ the Gaussian curvature, i.e., the product if the principal curvatures of $\partial D$ at the point $x.$ The Gaussian curvature coincides with the Jacobian of the Gaussian (spherical) map
$$\gamma: \partial D \to S^{n-1},$$
which maps any $x \in  \partial D$ to the external normal vector $\gamma(x)$ to $\partial D$ at $x.$
Respectively, the points of zero Gaussian curvature is exactly the set $Crit(\gamma)$ of critical points of the Gaussian map.

Pick any $x \in \partial D \setminus Crit(\gamma).$ Denote $\omega=\gamma(x).$
Since by Theorem \ref{T:convex} the boundary $\partial D$ is convex, all the principal curvatures $\kappa_1(x),...,\kappa_{n-1}(x)$ at $x$ are nonnegative. Since $ x \notin Crit(\gamma),$ the product $\kappa (x)=\kappa_1(x) \cdots \kappa_{n-1}(x) \neq 0$ and therefore $\kappa_1(x)>0,..., \kappa_{n-1}(x)>0.$

This implies that $V_D(\omega,t)$ has zero at $t=h_D^{+}(\omega)$ of order $\frac{n-1}{2}$ (see (\cite[1.7] {GG}) or the arguments of the proof
of Theorem \ref{T:even}). Therefore
$$V_D(\omega,t) =[h_D^{+}(\omega)-t]^{\frac{n-1}{2}} Q(t),$$
where $Q$ is a polynomial.

Repeating the argument for the point $x^{\prime} \in \partial D$ with the normal vector $\gamma(x^{\prime})=-\omega$ we conclude that
the polynomial $V_D$ has also zero of order $\frac{n-1}{2}$ at the point $t=h_D^{-}(\omega)=-h_D^{+}(-\omega).$ Therefore
\begin{equation}\label{E:zeros}
V_D(\omega,t) =[h_D^{+}(\omega)-t]^{\frac{n-1}{2}} [t-h_D^{-}(\omega]^{\frac{n-1}{2}} Q_1(t)
\end{equation}
and hence  $$\deg V_D(\omega, \cdot) \geq n-1.$$

By the choice of $x$, the estimate holds for all regular values $\omega=\gamma (x), x \in \partial D \setminus Crit(\gamma).$ By Sard 's theorem the Lebesgue measure $mes\{ \gamma(Crit(\gamma)\}=0$ and hence the estimate for the degree is fulfilled for almost all $\omega \in S^{n-1}.$
This proves the theorem.

\begin{theorem}\label{T:T1} Let $n$ be odd. Let $D$ be a  bounded domain in $\mathbb R^n$ with smooth boundary. Suppose that  for almost all $\omega \in S^{n-1}$ the function $V_D(\omega,t)$ is a polynomial in $t$, of degree at most $n-1.$ Then $D$ is an ellipsoid.
\end{theorem}

\pf
By our assumption and formula (\ref{E:zeros}),  for almost all $\omega \in S^{n-1}$ thd following representation holds:
$$V(\omega,t)= A(\omega) [(h_D^{+}(\omega)-t)(h_D^{-}(\omega)-t)]^{\frac{n-1}{2}}, h_D^{-}(\omega) < t < h_D^{+}(\omega),$$
where $A(\omega) >0.$

The function $V_D$ is the Radon transform of the characteristic function $\chi_D$ of the domain $D:$
$$V_D(\omega,t)= (R \chi_D)(\omega,t)=\int\limits_{<\omega,x>=t} \chi_D(x) dv_{n-1}(x)$$
Applying Projection-Slice Formula (\ref{E:SPF}) we obtain:
\begin{equation}\label{E:FT}
\begin{aligned}
\widehat \chi_D(r\omega)=\int\limits_{\mathbb R} e^{irt}V_D(\omega,t) dt=& \\
=&A(\omega) \int_{h_D^{-}(\omega)} ^{h_D^{+}(\omega)} e^{irt} [(h_D^{+}(\omega)-t)(h_D^{-}(\omega)-t)]^{\frac{n-1}{2}} dt.
\end{aligned}
\end{equation}

Consider the following functions on $S^{n-1}:$
$$B(\omega)=\frac{h_D^{+}(\omega)+h_D^{-}(\omega)}{2},$$
$$ C(\omega)=\frac{h_D^{+}(\omega)-h_D^{-}(\omega)}{2}.$$
The function $C(\omega)$ expresses the half-width of the domain $D$ in the direction $\omega.$

%%%%%

Performing  the following change of  variable in the integral in the left hand side in (\ref{E:FT}):
$$t=s+B(\omega),$$
we obtain
$$ \widehat \chi_D(r\omega)=(-1)^{\frac{n+1}{2}} A(\omega) e^{iB(\omega)r}  \int_{-C(\omega)} ^{C(\omega)} e^{irs} (C^2(\omega)-t^2)^{\frac{n-1}{2}} dt.$$
Performing the next  change of variable $ s=C(\omega)u$ we arrive at:
\begin{equation}\label{E:FT}
\widehat \chi_D(r\omega)=M(\omega)e^{iB(\omega)r} \int_{-1}^{1}e^{iC(\omega) ru}(1-u^2)^{\frac{n-1}{2}} du,
\end{equation}
where
$$ M(\omega)=(-1)^{\frac{n+1}{2}} A(\omega) (C(\omega))^{n-1}.$$

\begin{lemma} \label{L:MBC}
The function $M(\omega)$ is a (positive) constant. The function $B(\omega)$ is either zero or a homogeneous polynomial of degree one. The function $C^2(\omega)$ is  a homogeneous quadratic polynomial.
\end{lemma}

\pf of Lemma

The straightforward differentiation of both sides of (\ref{E:FT}) with respect to $r$ at the point $r=0$ gives:
\begin{equation}\label{E:MMBBCC}
\begin{aligned}
\widehat \chi_D(0)&=M(\omega), \\
\frac{d}{dr}\widehat \chi_D(r\omega)\vert_{r=0}&=M(\omega)iB(\omega)\alpha,\\
\frac{d^2}{dr^2}\widehat \chi_D(r\omega)\vert_{r=0}&=-M(\omega)\alpha[B^2(\omega)+C^2(\omega)],
\end{aligned}
\end{equation}
where we have denoted
$$\alpha:=\int_{-1}^{1}(1-u^2)^{\frac{n-1}{2}}du.$$

On the other hand,
\begin{equation}\label{E:mbc}
\begin{aligned}
\widehat \chi_D(0)&=vol(D), \\
\frac{d}{dr}\widehat \chi_D\vert_{r=0}&=i \int_D <\omega, y>dy, \\
\frac{d^2}{dr^2}\widehat \chi_D\vert_{r=0}&=-\int_D <\omega,y>^2dy.
\end{aligned}
\end{equation}

Therefore, from (\ref{E:MMBBCC}) and the first equality in (\ref{E:mbc}) we obtain
$$
M(\omega)=vol (D)>0.
$$
By (\ref{E:MMBBCC}) and the second equality in (\ref{E:mbc})  we have
$$
B(\omega)=\frac{1}{\alpha vol(D)} \int_D <\omega, y>dy,
$$
and hence $B$ is a linear form. In particular, it can be identically zero, if, for instance, the domain  $D$ has a  central symmetry.

Finally, (\ref{E:MMBBCC}), the third equality in (\ref{E:mbc}) imply
$$B^2(\omega)+C^2(\omega)=\frac{1}{\alpha vol(D)} \int_D <\omega,y>^2 dy$$ and then
$$C^2(\omega)=\frac{1}{(vol D)^2} [\int\limits_D <\omega,y>^2 dV(y) - B^2(\omega)$$
is a quadratic form, because $B$ is a linear form.  The lemma is proved.

Let us proceed with the proof of Theorem \ref{T:T1}. Applying a suitable orthogonal transformation we can assume that the quadratic form $C^2(\omega)$ has a diagonal form:
$$C^2(\omega)=c_1\omega_1^2 +...+c_n \omega_n^2.$$
Also, according to Lemma \ref{L:MBC} one has
$$M(\omega)=M=const \ \mbox{and} \  B(\omega)=<b,\omega>,$$
for some vector $b.$

Remember that
$$
\begin{aligned}
&B(\omega)=\frac{1}{2}(h_{+}(\omega)+h_-(\omega)), \\
& C(\omega)=\frac{1}{2}(h_+(\omega)-h_{-}(\omega)).
\end{aligned}
$$
Applying translation by vector $b$ and a dilatation of the domain $D,$ and using property (\ref{E:D+a}),  we may assume that $M=1, B(\omega)=0.$
This means that $h_D^{+}=h_D^{-}$ and then $C=h_D^+.$
Thus,
$$h_+(\omega)=C(\omega)=\sqrt{\sum_{j=1}^n c_j \omega_j^2}.$$
But the function in the right hand side  is exactly the support function $h_E^{+}$ of the ellipsoid
$$E=\{\sum_{j=1}^n \frac{x_j^2}{c_j}=1\}$$ and hence
$$D=E.$$ This completes the proof of Theorem \ref{T:T1}.

The following result is a combination of Theorem \ref{T:T1} and Lemma \ref{L:a_k}:
\begin{theorem}\label{T:T2} Let $D$ be a domain in $\mathbb R^n, $  with a smooth boundary, where $n$ is odd. Suppose that the section-volume function $V_D$ is a polynomial,
$V_D(\omega,t)=a_0(\omega)+a_1(\omega)t + \cdots a_N(\omega)t^N, \ h_D^{-}(\omega) \leq t \leq h_D^{+}(\omega)$
and for each $k > n-1$ the coefficient $a_k(\omega)$  coincides with the restriction to $S^{n-1}$ of a polynomial of degree at most $k-n+1.$ Then $D$ is an ellipsoid.
\end{theorem}
\pf
Lemma \ref{L:a_k} and the condition imply that $a_k=0$ for $k>n-1.$ Theorem \ref{T:T1} implies that $D$ is an ellipsoid.

\section {Axially symmetric polynomially integrable domains in 3D are ellipsoids}

In this section we confirm the conjecture about ellipsoids as the only polynomially integrable domains, however under the pretty strong condition of axial symmetry.
To avoid additional technical difficulties, we consider the three-dimensional case.

The following theorem asserts that, in this case,  to conclude that the domain is an ellipsoid,
it suffices essentially to demand the polynomial integrability only in two directions:

\begin{theorem} \label{T:axial} Let the $D$ be an axially and centrally symmetric  domain in $\mathbb R^3.$  Suppose that
the section-volume function $V_D(\omega,t)$ is a polynomial in $t,$ when $\omega$ is one of the two orthogonal directions $\omega=\xi$ and $ \omega=\eta,$
where $\xi$ is the symmetry axis  and $\eta$ is any orthogonal direction.
Assume also that $V_D(\xi,t)=b_0 + b_1t+ \cdots +b_{2N}t^{2N}$ and $ b_{2N} <0.$
Then $D$ is an ellipsoid.
\end{theorem}

\pf
We assume that the center of symmetry of $D$ is 0 and the symmetry axis $\xi$ the $z-$ axis, $\xi=(0,0,1).$

Due to the axial symmetry, the cross-section of $D$ by the plane $z=t$ is a two-dimensional disc
$D(0, r(t))$ of radius $r(t)$ and therefore the volume of the corresponding cross-section is
$$V_D(\xi,t)= \pi[r(t)]^{2}.$$
.
By our assumption, $V_D(\xi,t)$ is a polynomial
$$P(t):=V_D(\xi,t).$$
The central symmetry implies that it is an even polynomial, i.e.
$$V_D(\xi,t)=P(t)=b_0 + b_2t^2+ \cdots +b_{2N}t^{2N}$$
and hence
$$P(t)= \pi [r(t)]^2.$$

Since $t=z$ and $r(t)^2=x^2+y^2$, we conclude that the equation of the boundary $\partial D$ of the domain $D$ is given by:
\begin{equation}\label{E:eq}
x^2 + y^2 =P(z)=b_0 + b_2z^2+ \cdots +b_{2N}z^{2N},
\end{equation}
where we have incorporated the constant $\frac{1}{\pi}$ into $P.$

Consider now the cross-section of the domain $D$ in the orthogonal direction $\eta.$ This vector lies in $x,y-$ plane and due to rotational symmetry can be taken $\eta=(0,1,0).$
Intersect $D$ by the hyperplane $y=t.$ The equation of that section can be written as:
\begin{equation} \label{E:section}
x^2 -(b_2 z^2+...+b_{2N} z^{2N} )= b_0-t^2.
\end{equation}
Now introduce the family $\Omega(\alpha), \alpha >0$ of domains in $\mathbb R^2_{x,z}:$
$$\Omega(\alpha)=\{(x,z) \in \mathbb R^2: x^2 -(b_2 z^2+...+b_{2N} z^{2N} ) \leq \alpha^2 \} .$$
Every domain $\Omega(\alpha)$ is non-empty, as $0 \in \Omega(\alpha).$ Also, the domains $\Omega(\alpha)$ are bounded
because the condition $b_{2N}<0$ implies that the left hand side of the inequality tends to $\infty$ when $(x,z) \to \infty$ and hence
$\mathbb R^2 \setminus B_R$ is free of points from $\Omega(\alpha),$ for any fixed $\alpha$ and $R$ large enough.

Note that $b_0= vol(D \cap \{z=0\}) >0.$  The parameter $t=y$ varies in a neighborhood of $0$
and therefore $b_0-t^2 >0$ for small $t$ and, correspondingly,
$$D \cap \{y=t\}=\Omega(\sqrt{b_0-t^2}),$$
for $b_0 -t^2>0.$
The function
$ A(\alpha):=vol \Omega(\alpha)$
is real analytic in $\alpha.$ On the other hand,due to the central symmetry the function
$$vol (D \cap \{y=t\})=V_D(\eta, t)$$
is a polynomial, $B(t^2),$ in $t^2,$

Since
$$A(\alpha)=B(b-\alpha^2)$$
for $\alpha >0$ from an open set, we conclude that the volume function $A(\alpha)=vol \Omega(\alpha)$ is a polynomial on the whole half-line $\alpha >0.$

On the other hand, by the change of variables
$$u=\frac{x}{\alpha}, v=\frac{z}{\alpha^{\frac{1}{N}}}$$
we obtain
$$
\begin{aligned}
vol \Omega(\alpha)=\int\limits_{ x^2 -(b_2 z^2+...+b_{2N} z^{2N} ) \leq \alpha^2 }dxdz &= \\
&=  \alpha ^{1+\frac{1}{N}} \int\limits_{u^2-(b_2v^2 \alpha^{\frac{2}{N}-2}+...+b_{2N}v^{2N}) \leq 1 } du dv.
\end{aligned}
$$
Since
$$ \lim_{\alpha \to \infty}  \int\limits_{u^2-(b_2v^2 \alpha^{\frac{2}{N}-2}+...+b_{2N}v^{2N}) \leq 1 } du dv=
\int\limits_{u^2- b_{2N}v^{2N} \leq 1 } du dv < \infty,
$$
we conclude that
$$A(\alpha)=vol \Omega(\alpha)= O(\alpha^{1+\frac{1}{N}}), \alpha \to \infty.$$
Since $A(\alpha)$ is a polynomial, the exponent $1+\frac{1}{N}$ must be integer, which implies $ N=1.$

Thus, the equation (\ref{E:eq}) of the domain $D$ is given by the second order equation
$$ x^2 + y^2 = b_0 + b_2z^2, b_2 <0,$$
and $D$ is an ellipsoid.

\section{Finite stationary phase expansion point of view}\label{S:Stat}

Let $D$ be a convex polynomially integrable domain with the section-volume function
$$V_D(\omega,t)=\sum\limits_{j=0}^N a_j(\omega)t^j.$$

The relation (\ref{E:SPF}) between Fourier and Radon transforms yields
\begin{equation}\label{E:for}
\begin{aligned}
\widehat \chi_D(r\omega)= \sum\limits_{j=0}^N a_j(\omega)\int\limits_{h_D^{-}(\omega)}^{h_+(\omega)}t^j e^{irt}dt=& \\
&=\sum\limits_{j=0}^N a_j(\omega)\frac{d^j}{dr^j}\frac{ e^{h_D^+(\omega) r}-e^{h_D^{+}(\omega)r }}{r}.
\end{aligned}
\end{equation}
On the other hand, by the Green's formula
\begin{equation}
\begin{aligned}
\widehat \chi_D (r\omega)&=\int\limits_D e^{ir<\omega,x>}dV(x)=& \\
=& - \frac{1}{r^2}\int\limits_D \Delta e^{ir<\omega,x>}dV(x)=\\
=& - \frac{1}{r^2}\int\limits_{\partial D} e^{ir<\omega,x>}<\omega,\nu(x)>dS(x),
\end{aligned}
\end{equation}
where $dS$ is the surface measure on $\partial D$ and $\nu(x)$ is the external unit normal vector to $\partial D.$

The right hand side in (\ref{E:for}  can be written, after performing differentiation in $r,$ in the form
$$e^{irh_D^{+}(\omega)} Q^+(\omega, \frac{1}{r})+  e^{irh_D^{-}(\omega)} Q^{-}(\omega, \frac{1}{r}),$$
where $Q^{\pm}(\omega,t)$ are polynomials with respect to $t.$
Thus, we have
\begin{equation}\label{E:stat}
\int\limits_{\partial D} e^{ir<\omega,x>}<\omega,\nu(x)>dS(x)=e^{irh_D^{+}(\omega)} r^2Q^+(\omega, \frac{1}{r})-  e^{irh_D^{-}(\omega)} r^2 Q^{-}(\omega, \frac{1}{r}).
\end{equation}

The critical points of the function
$$ \varphi_{\omega}(x):= <\omega,x>$$
are the points $x^+$ and $x^-$ such that
\begin{equation}
\begin{aligned}
<\omega,x^+> &= \max\limits_{x \in \overline D}<\omega,x>=h_{D}^+(\omega) ,  \\
<\omega,x^-> &= \min\limits_{x \in \overline D}<\omega,x>=h_{D}^{-}(\omega).
\end{aligned}
\end{equation}
Equation (\ref{E:stat}) is the stationary phase expansion (cf. \cite{Wong}) in the powers of $\frac{1}{r}$
of the oscillatory integral in the left hand side, over the manifold $\partial D$ with the large parameter $r,$ at the critical points $x^{\pm}$ of the phase function $\varphi_{\omega}(x).$

In a general setting, the stationary phase expansion is an asymptotic series. In our particular case, however, there is only finite number of  terms in the expansion, which is rather an exact equality.

The phenomenon of stationary phase expansions with finite number of terms was studied in \cite{Bott},\cite{Duist}, \cite{Bernard}.
The simplest example of such expansion delivers the unit ball $B^3$  in $\mathbb R^3.$ The Fourier transform of the characteristic function of the ball is
\begin{equation}
\begin{aligned}
\frac{1}{\pi}\widehat \chi_{B^n}(r\omega)&= \int\limits_{-1}^{1} e^{irz}(1-z^2)dz=&  \\
=& e^{ir}(\frac{1}{ir}-\frac{2}{(ir)^2}+\frac{2}{(ir)^3}) + e^{-ir}(-\frac{1}{ir}-\frac{2}{(ir)^2}-\frac{2}{(ir)^3}).
\end{aligned}
\end{equation}
The right hand side is the stationary phase expansion, at the critical points $z=\pm 1,$ of the oscillatory integral on $S^2,$ resulted from the Green's formula for the integral over the ball in the left hand side. The phase function of this oscillatory integral is the height-function $\varphi (x,y,z)=z.$

It was shown in \cite{Duist} that finite (single-term) stationary phase expansion appears for the phase functions - moment maps of Hamiltonian actions on symplectic manifolds.
This result was generalized in \cite{Bernard} for the even-dimensional, not necessarily symplectic, manifolds acted on by $S^1$ with isolated fixed points. The situation studied in \cite{Bernard} leads to multi-term expansions.

It is an open problem to describe all the situations when the stationary phase expansions are finite (e.g., see the discussion in \cite{Bernard}). Balls and ellipsoids do deliver such examples, but apparently do not exhaust them. We saw that polynomially integrable domains and the phase functions \ \ $\varphi_{\omega}(x)=<\omega, x>$ (the height-function in the direction $\omega$) also provide finite stationary phase expansions (\ref{E:stat}).

However, in the case of polynomially integrable domains the situation is much more restrictive. Indeed, we have finite expansions, on the  manifold $\partial D,$ but for the whole family $\omega-$ height functions $\varphi_{\omega}, $ parametrized by the unit vector $\omega \in S^{n-1}.$

In fact, this property fully characterizes   polynomially integrable domains:
\begin{proposition} Let $D$ be a convex domain in $\mathbb R^n$ with smooth boundary. Then $D$ is polynomially integrable if and only if the Fourier transform
$\widehat \chi_D(r\omega)$ has, for any direction $\omega \in S^{n-1},$ a  finite stationary phase expansion with respect to the large parameter $r.$
\end{proposition}

\pf

The part "if" is already proved by formula (\ref{E:stat}). We have to show that if (\ref{E:stat}) is valid then $D$ is a polynomially integrable domain.
By Projection-Slice Formula (\ref{E:SPF}) the section-volume function $V_D(\omega,t)$ is the inverse one-dimensional Fourier transform with respect to $r$ of $\widehat \chi(r\omega):$
$$
V_D(\omega,t)= F^{-1}_{r \to t}\{  e^{ih_D^+(\omega) r}Q^{+}(\omega,\frac{1}{r})+e^{ih_D^{-}(\omega) r}Q^{-}(\omega,\frac{1}{r}) \}.
$$
The general terms in the right hand side are
$$ v_k^{\pm}(\omega,r):=q^{\pm}_k(\omega) e^{ih_D^{\pm}(\omega) r}\frac{1}{r^k}.$$
Consider one of them, say, $v_k^+.$

If $$f_k(\omega, t)=F^{-1}_{r \to t} v_k^{+}(\omega, r)$$
then its $k$th derivative in $t$ is
\begin{equation}
\begin{aligned}
\frac{d^k}{dt^k} f_k (\omega,t)=F^{-1}_{r \to t} (ir)^k v_k(\omega,r)&= \\
&=i^k q^{+}_k(\omega)F^{-1}_{r \to t} e^{irh_D^{\pm}(\omega) }=i^k \delta(t-h^{+}(\omega)).
\end{aligned}
\end{equation}
Consequent integration in $t$ yields that $f_k(\omega,t)$ is a piecewise polynomial, namely, it is a polynomial in $t$
on  $t<h_D^{+}(\omega),$ and is a polynomial in $t$ on $t > h_D^{+}(\omega).$

Likewise, the functions $F^{-1}_{r \to t}v_k^{-}(\omega,t)$ are polynomials on each half-line $t<h_D^{-}(\omega)$ and $t > h_D^{-}(\omega).$
Then the section-volume function $V_D$ is  a polynomial on the segment $h_D^{-}(\omega) \leq t \leq h_D^{+}(\omega),$  as a sum of $F^{-1}v_k^{\pm}.$ Thus, $D$ is polynomially integrable.

\section{Open questions}\label{S:Open}
 We conclude with formulating  open questions.
\begin{itemize}
\item Are ellipsoids  the only domains with infinitely smooth boundaries in the Euclidean spaces $\mathbb R^n$ of odd dimensions $n$ for which
the section-volume function $V_D(\omega,t)$ is a polynomial in $t$? The same question refers to the cut-off function $\widehat V(\omega,t)$ (see \cite{Vas}).
\item Under which conditions for the behaviour of the above volume functions with respect to $\omega$   polynomially integrable domains are ellipsoids?
\item What geometric properties of a domain can be derived from the finite stationary phase expansion for Fourier transform of its characteristic function?
\item What is the role of the smoothness condition for the boundary? Are there  polynomially integrable domains with a relaxed condition of smoothness, different from ellipsoids?
\end{itemize}

\bigskip
%\section{Acknowledgements}
I am grateful to Mikhail Zaidenberg for drawing my attention to the subject of this article and stimulating initial discussions.

Bar-Ilan University, Holon Institute of Technology; Israel

{\it E-mail}: agranovs@math.biu.ac.il    

\end{document}